\documentclass[10pt, titlepage]{amsart}

\usepackage{ae} 

\usepackage[cp1250]{inputenc}
\usepackage{amsmath}
\usepackage{amssymb, amsfonts,amscd,verbatim}

\usepackage{latexsym}
\usepackage{tikz}
\usepackage{indentfirst}
\usepackage{latexsym}

\usepackage{verbatim}   
\usepackage{amsmath}    

\theoremstyle{plain}
\newtheorem{Prop}{Proposition}[section]
\newtheorem{Thm}[Prop]{Theorem}
\newtheorem{Cor}[Prop]{Corollary}

\newtheorem{Lem}[Prop]{Lemma}

\theoremstyle{definition}
\newtheorem{Def}[Prop]{Definition}

\theoremstyle{remark}
\newtheorem{Rem}[Prop]{Remark}

\newtheorem{Example}[Prop]{\bf Example}

\def\int{\mathop{\roman{int}}}

\def\1{^{-1}}

\def\RR{{\mathbf R}}
\def\NN{{\mathbf N}}

\def\im{\text{Im}}
\def\dim{\text{dim}}

\def\diam{\text{diam}}

\def\asdim{\mathrm{asdim}}

\def\dim{\mathrm{dim}}
\def\diam{\mathrm{diam}}
\def\dokaz{{\bf Proof. }}
\def\edokaz{\hfill $\blacksquare$}

\def\a{\alpha}
\def\b{\beta}
\def\g{\gamma}
\def\d{\delta}

\def\f{\varphi}

\def\z{\zeta}

\def\l{\lambda}
\def\m{\mu}

\def\di{\partial}

\numberwithin{equation}{section}
\begin{document}
\title[Inducing maps between Gromov boundaries]{Inducing maps between Gromov boundaries}%

\author{J.~Dydak}
\address{University of Tennessee, Knoxville, TN 37996, USA}
\email{jdydak@vol.utk.edu}

\author{\v Z.~Virk}
\address{University of Ljubljana, Jadranska 19, SI-1111 Ljubljana, Slovenia}
\email{zigavirk@gmail.com}

\date{ \today
}

\thanks{This research was partially supported by the Slovenian Research
Agency grants P1-0292-0101.}

\begin{abstract}
It is well known that quasi-isometric embeddings of Gromov hyperbolic spaces induce
topological embeddings of their Gromov boundaries. A more general question is to detect
classes of functions between Gromov hyperbolic spaces that induce continuous
maps between their Gromov boundaries. In this paper we introduce the class of visual
functions $f$ that do induce continuous
maps $\tilde f$ between Gromov boundaries. Its subclass, the class of radial functions,
induces H\"older maps between Gromov boundaries. Conversely, every H\"older map
between Gromov boundaries of visual hyperbolic spaces induces a radial function.
We study the relationship between large scale properties of $f$
and small scale properties of $\tilde f$, especially related to the dimension theory.
In particular, we prove a form of the dimension raising theorem. We give a natural example of a radial dimension raising map
and we also give a general class of radial functions that raise asymptotic dimension.
\end{abstract}

\maketitle

\section{Introduction}

Gromov boundary (as defined by Gromov) is one of the central objects in the geometric group theory and plays a crucial role in the Cannon's conjecture \cite{C}. It is a compact metric space that represents an image at infinity of a hyperbolic metric space. It is well known that quasi-isometries (and quasi-isometric embeddings) induce homeomorphisms (respectively, embeddings) of Gromov boundaries.  If a continuous map between Gromov boundaries is induced by a proper embedding it is called a  Cannon-Thurston map \cite{Mi}.  Baker and Riley \cite{BR} showed that a proper embedding may not induce a Cannon-Thurston map. \' Alvaro Mart\' inez-P\' erez \cite{MP}  proves that a power quasi-symmetric(or PQ-symmetric) homeomorphism between two complete metric spaces can be extended to a quasi-isometry between their hyperbolic approximations.
\' Alvaro Mart\' inez-P\' erez and Manuel A. Mor\' on \cite{MM1} discuss uniformly continuous maps between ends of R-trees (see also
\cite{MM2}).
Bruce Hughes, \' Alvaro Mart\' inez-P\' erez and Manuel A. Mor\' on \cite{HMM} have results on bounded distortion homeomorphisms on ultra metric spaces.
In this paper we introduce radial functions that represent a subclass of large-scale Lipschitz functions inducing continuous maps on the Gromov boundary. We prove that a number of coarse properties of functions transfer via radial functions into their counterparts, the classical properties of continuous maps. Such properties include (coarse) surjectivity and (coarse) n-to-1 property. As a consequence, we present an improved dimension raising theorem for radial functions (see \cite{MV} for the original theorem).

\begin{Thm}[Large scale dimension raising theorem for hyperbolic spaces]
Suppose $X$, $Y$ are proper $\d-$hyperbolic geodesic spaces and $f\colon X \to Y$ is a radial function. If $f$ is coarsely $(n+1)$-to-1 and coarsely surjective then
$$
\dim(\di Y)\leq \dim (\di X)+n.
$$

Moreover, if $X$ and $Y$ are hyperbolic groups then
$$
\asdim(Y)\leq \asdim (X)+n.
$$
\end{Thm}

In the penultimate section we provide a natural example of a situation as described in the main theorem. It is essentially a nice embedding of the Cayley graph of the free group on two generators into the hyperbolic plane.

Coarsely n-to-1 maps were introduced in \cite{MV}. They appear naturally as quotient maps via a finite group action by coarse functions on a metric space. Quotient maps via a group action are usually considered in the context of the \v Svarc-Milnor Lemma where the action is required to be cocompact. The study of coarsely n-to-1 functions is the only generalization to a non-cocompact action for which meaningful connections between the acting group and the space acted upon has been made.

In the last section of the paper we investigate a general method of constructing coarsely $n$-to-1 functions of proper hyperbolic spaces
via extensions of $n$-to-1 maps of their boundaries.

\section{Sequential boundaries and visual maps}
We recall a number of well established concepts. For equivalent definitions of a hyperbolic space see \cite{BS}. Beware that various definitions of hyperbolicity are equivalent even though the parameter $\delta$ may change by a multiplicative constant. In order to make use of a definition where parameter $\delta$ fits into all definitions of hyperbolicity we define (see Definition \ref{DefBasic}) $\delta$-hyperbolicity by a $\delta/4$-inequality instead of a $\delta$-inequality (see Proposition 2.1.3 of \cite{BS} for one such reason).
\begin{Def}\label{DefBasic}
Given a metric space $(X,d)$, the \textbf{Gromov product} of $x$ and $y$ with respect to $a\in X$ is defined by
$$
(x,y)_a=\frac{1}{2}\big(d(x,a)+d(y,a)-d(x,y)\big).
$$

A metric space $(X,d)$ is (Gromov) $ \delta-$\textbf{hyperbolic} if it satisfies the $\delta/4$-inequality:
$$
(x,y)_{x_0} \geq \min \{(x,z)_{x_0},(z,y)_{x_0}\}-\delta/4, \quad \forall x,y,z,x_0\in X.
$$

A metric space is \textbf{proper} if closed balls are compact.

Given points $x,y$ in a geodesic space, a \textbf{geodesic} between $x$ and $y$ is denoted by $xy$. We refrain from distinguishing the geodesic as a map and the image of the geodesic when the context of usage is apparent.
\end{Def}

\begin{Def}\label{DefSequentialBoundary}\cite{KB}
Suppose $X$ is a proper, hyperbolic geodesic space and $a\in X$ is a basepoint.

The \textbf{sequential boundary} $\partial X$ consists of classes of sequences going "straight" to infinity:
$$
\di X=\{(x_n) \mid (x_n) \emph{ sequence in } X, \liminf_{i,j\to \infty} (x_i,x_j)_a=\infty \} / \sim
$$
where the equivalence of sequences is given by
$$
((x_n) \sim (y_n)) \Leftrightarrow (\liminf_{i,j\to \infty} (x_i,y_j)_a=\infty).
$$
The topology on $\partial X$ is induced by sets $\{U(p,r)\mid r>0, p\in \di X\}$, where
$$
U(p,r)=\{q\in \di X \mid \exists \emph{ sequences } (x_n), (y_n):
$$
$$
[(x_n)]=p, [(y_n)]=q, \liminf_{i,j\to \infty} (x_i, y_j)_a \geq r\}.
$$
\end{Def}

\begin{Rem} We can extend the sets $U(p,r)$ to be subsets of $X\cup\partial X$ as follows:
$x\in U(p,r)$ if and only if for $[(x_n)]=p$ one has $\liminf_{i\to \infty} (x_i, x)_a \geq r$.
This way one can create a compact topological space $X\cup\partial X$.
\end{Rem}

\begin{Def}
Suppose $X$ and $Y$ are metric spaces. A function $f\colon X \to Y$ is \textbf{large scale Lipschitz} (LSL, or $(\l_1,\m_1)-$LSL) if there exist $\l_1,\m_1>0$ so that

$$
 d(f(x),f(y))\leq \l_1 \cdot d(x,y)+\m_1, \qquad \forall  x,y\in X.
$$
\end{Def}

\begin{Def}\label{DefVisualMap}
A large scale Lipschitz function $f:X\to Y$ of Gromov hyperbolic spaces is called \textbf{visual} if
$$(x_n,y_n)_a\to\infty \implies (f(x_n),f(y_n))_{f(a)}\to\infty$$
for some $a\in X$.
\end{Def}

\begin{Rem}
 Since $(x,y)_a-2d(a,b)\leq (x,y)_b\leq (x,y)_a+2d(a,b)$ for all $a,b,x,y\in X$,
 $(x_n,y_n)_a\to\infty \implies (f(x_n),f(y_n))_{f(a)}\to\infty$ is equivalent to
 $(x_n,y_n)_b\to\infty \implies (f(x_n),f(y_n))_{f(b)}\to\infty$ and the definition of $f$ being visual
 does not depend on the base point $a$.
\end{Rem}

\begin{Prop}\label{BasicVisualProp}
Suppose $X,Y$ are Gromov hyperbolic spaces, $f:X\to Y$ is a large scale Lipschitz function,
and $f(a)=b$. $f$ is visual if
 and only if for each $r > 0$ there is $s > 0$ such that $(x,y)_a > s$ implies $(f(x),f(y))_{b} > r$.
\end{Prop}
\dokaz By contradiction, assume $r >0$ and there is a sequence $(x_n,y_n)\in X\times X$
such that $(x_n,y_n)_a > n$ but $(f(x_n),f(y_n))_{b} \leq r$. However, in this case
$(x_n,y_n)_a\to\infty$ and $(f(x_n),f(y_n))_{f(a)}$ is bounded, a contradiction.
\edokaz

\begin{Thm} [Induced continuous maps on the boundary]
Suppose $f:X\to Y$ is a large scale Lipschitz function of proper hyperbolic spaces. If $f \colon X \to Y$ is visual, then $f$ induces a function
      $$
      \tilde f\colon X\cup\partial X \to Y\cup \partial Y
      $$
      such that $\tilde f(\partial X)\subset \partial Y$ and $\tilde f$ is continuous at all points of $\partial X$.

      Conversely, if $f$ induces a function
      $$
      \tilde f\colon X\cup\partial X \to Y\cup \partial Y
      $$
      such that $\tilde f(\partial X)\subset \partial Y$ and $\tilde f$ is continuous at all points of $\partial X$, then $f$ is visual.
\end{Thm}
\dokaz
Fix $a\in X$. Consider the sequential boundary of Definition \ref{DefSequentialBoundary}. Given a sequence $(x_n)$ define $\tilde f (x_n)=(f(x_n))$. The map $\tilde f$ is well defined:
\begin{itemize}
  \item if $\liminf _{i,j\to\infty}(x_i,x_j)_{a}=\infty$ then
    $$
    \liminf _{i,j\to\infty}(f(x_i),f(x_j))_{f(a)}=\infty
    $$
Hence $\tilde f$ is well defined on every single sequence.
  \item if $\liminf _{i,j\to\infty}(x_i,y_j)_{a}=\infty$ then
    $$
    \liminf _{i,j\to\infty}(f(x_i),f(y_j))_{f(a)}=\infty
    $$
Hence $\tilde f$ is well defined on the equivalence classes of sequences, i.e., on the boundary.
\end{itemize}

To prove $\tilde f$ is continuous at points of $\partial X$ use the fact that for
each $r > 0$ there is $s > 0$ such that
$(x,y)_a > s$ implies $(f(x),f(y))_{f(a)} > r$ (see \ref{BasicVisualProp}).
Therefore we have $\tilde f (U(p,s))\subset U(\tilde f(p),r), \forall p\in \partial X$. (Notice that sets $U(p,r)$ are decreasing as $r$ is increasing.)

Assume $f$ induces a function
      $$
      \tilde f\colon X\cup\partial X \to Y\cup \partial Y
      $$
      such that $\tilde f(\partial X)\subset \partial Y$ and $\tilde f$ is continuous at all points of $\partial X$.
      Fix $a\in X$ and assume existence of $x_n,y_n\in X$ such that $(x_n,y_n)_a\to \infty$ but $(f(x_n),f(y_n)_{f(a)} \leq M < \infty$
      for all $n$. Due to $X$ being proper, we can reduce the argument (by switching to subsequences) to the case
      of $p=[(x_n)]$ and $q=[(y_n)]$ belonging to $\partial X$.
      Now, one has $p=\lim x_n=\lim y_n=q$ in $X\cup\partial X$. Therefore $\tilde f(p)=\lim f(x_n)$ in $Y\cup\partial Y$
      and $\tilde f(q)=\lim f(y_n)$ in $Y\cup\partial Y$ which implies $(f(x_n),f(y_n)_{f(a)}\to \infty$, a contradiction.
      Thus, $f$ is visual.
\edokaz

\begin{Rem}
  A similar result for proper embeddings was proved in \cite[Lemma 2.1]{Mi} .
\end{Rem}

\section{Geodesic boundaries and radial functions}

In this section we recall the concept of the geodesic boundary of a hyperbolic space
and we introduce geometrically the concept of a radial function between hyperbolic spaces.
We prove that radial functions are visual and induce H\" older maps on geodesic
boundaries (when equipped with visual metrics).

\begin{Def}\cite{KB}
Suppose $X$ is a proper, hyperbolic geodesic space.

The geodesic boundary consists of classes of geodesic rays:
$$
\di^g_x X=\{\gamma \mid \gamma \colon [0,\infty)\to X, \emph{ geodesic }, \gamma(0)=x\} / \sim
$$
where the equivalence of rays is given by
$$
(\gamma_1 \sim \gamma_2) \Leftrightarrow (\exists K: d(\g_1(t),\g_2(t))<K, \forall t).
$$
The topology consists of basis $\{V(p,r)\mid r>0, p\in \di^g_x X\}$, where
$$
V(p,r)=\{q\in \di^g_x X \mid \exists \emph{ geodesic rays } \g_1, \g_2 \emph{ based at x: }
$$
$$
[\g_1]=p, [\g_2]=q, \lim_{t\to \infty} (\g_1(t),\g_2(t))_x \geq r\}.
$$

\end{Def}

It turns out that all definitions are independent of $x$. Furthermore, all boundaries of a fixed proper $\d$-hyperbolic geodesic space $X$ are naturally homeomorphic to each other as compact metric spaces \cite[Proposition 2.14]{KB}. Throughout the paper we will thus frequently (and without a particular notice) use geodesic and sequential representation of points on the boundary and switch between them.

The following definition introduces a new type of functions called the radial functions. They represent a subclass of radial functions as will be shown later.

\begin{Def}\label{BasicIsoDefs}
Suppose $X$ and $Y$ are metric spaces. A function $f\colon X \to Y$ is a \textbf{quasi-isometric embedding} if there exist $\l,\m>0$ so that

$$
\l^{-1} \cdot d(x,y)-\m \leq d(f(x),f(y))\leq \l \cdot d(x,y)+\m, \qquad \forall  x,y\in X.
$$

 $f\colon X \to Y$ is a \textbf{quasi-isometry} if it is a quasi-isometric embedding and $f(X)$ is coarsely dense in $Y$, i.e., $\exists R: \forall y\in Y \exists x\in x: d(y,f(x))\leq R$.

 $f\colon X \to Y$ is \textbf{radial} ($(\l_2,\m_2)-$radial) with respect to $x_0\in X$ if it is large scale Lipschitz and there exist $\l_2,\m_2>0$ so that for every finite geodesic $\g\colon [0,M]\to X, \g(0)=x_0$ the following condition holds:
$$
\l_2 \cdot d(x,y)-\m_2 \leq d(f(x),f(y)), \qquad \forall  x,y\in \im \g.
$$
\end{Def}

Note that every quasi-isometry is a radial function. The converse does not hold as $\RR\to \RR; x\mapsto |x|$ demonstrates. However, it is apparent that every radial LSL function $[0,\infty)\to X$ is a quasi-isometric embedding.

\begin{Prop}  \label{PropIndependenceOfQRofBasepoint}
Suppose $f\colon X\to Y$ is a $(\l_1,\m_1)-$LSL function and $X$ is a $\d-$hyperbolic space.  If $f$ is radial with respect to some point, then it is radial with respect to any other point.
\end{Prop}

\dokaz
Suppose $f$ is $(\l_2,\m_2)-$radial with respect to $x_0.$ We will prove it is also radial with respect to $x_1$.

Choose $x,y$ on some geodesic starting at $x_1$ so that $x\in x_1 y$. There exists a point $\tilde x$ on the union of geodesics $x_0 y \cup x_1 x_0$ so that $d(x,\tilde x)\leq \d$.

Case 1: $\tilde x \in x_0 y$. Note that $d(f(x),f(\tilde x))\leq \l_1 d(x,\tilde x)+\m_1\leq \l_1 \d + \m_1$ hence
$$
d(f(x),f(y)) \geq d(f(y),f(\tilde x))-d(f(\tilde x),f(x))\geq
$$
$$
\geq  \l_2 \cdot d(y,\tilde x) -\m_2 -\l_1 \d - \m_1 \geq \l_2 \cdot (d(y, x)-d(x,\tilde x)) -(\m_2 +\l_1 \d +\m_1)\geq
$$
$$
\geq \l_2 \cdot d(y, x) -(\l_2 \d +\m_2 +\l_1 \d +\m_1).
$$

Case 2: $\tilde x \in x_0 x_1$. Note that $d(x,x_0)\leq \d + d(x_0,x_1)$ and $d(f(x),f(x_0))\leq \l_1(\d + d(x_0,x_1))+\m_1$ hence
$$
d(f(x),f(y)) \geq d(f(y),f( x_0))-d(f( x_0),f(x))\geq
$$
$$
\geq \l_2 d(y,x_0)-\m_2-\l_1(\d + d(x_0,x_1))-\m_1\geq
$$
$$
\geq\l_2 (d(y,x)-d(x,x_0))-\m_2-\l_1(\d + d(x_0,x_1))-\m_1=
$$
$$
=\l_2 (d(y,x)-d(x,x_0))-(\m_2+\l_1(\d + d(x_0,x_1))+\m_1)=
$$
$$
= \l_2 d(y,x)-(\l_2d(x,x_0)+\m_2+\l_1(\d + d(x_0,x_1))+\m_1)\geq
$$
$$
\geq \l_2 d(y,x)-(\l_2(\d +d(\tilde x,x_0))+\m_2+\l_1(\d + d(x_0,x_1))+\m_1)\geq
$$
$$
\geq \l_2 d(y,x)-(\l_2(\d +d(x_1,x_0))+\m_2+\l_1(\d + d(x_0,x_1))+\m_1).
$$
Altogether we have proved that $f$ is a radial function with respect to $x_1$ as the obtained parameters are independent of $x$ and $y$.
\edokaz

\begin{Example} Notice it is insufficient to define radial functions using infinite rays only instead of finite (bounded) geodesics.
Indeed, consider two different rays $r_1$ and $r_2$ stemming from a point in $H^2$. Create a tree $T$ consisting of one infinite
ray $r$ and countably many geodesics $g_n$ joining $r(n)$ to a point $x_n$ at the distance equal to $dist(r_1(n),r_2(n))$.
 The map $f:T\to H^2$ sends $r$ isometrically onto $r_1$ and each $g_n$ isometrically onto the geodesic from
 $r_1(n)$ to $r_2(n)$. It is clear that $f$ is $1$-Lipschitz, $T$ has only one boundary point,
 and $f$ is not visual despite the fact that for every infinite geodesic $\g\colon [0,\infty)\to T, \g(0)=x_0$ the following condition holds:
$$
1\cdot d(x,y)-0 \leq d(f(x),f(y)), \qquad \forall  x,y\in \im \g.
$$

\end{Example}

However, in case of visual hyperbolic spaces we can restrict our attention to geodesic rays.
\begin{Def}
 A hyperbolic space $X$ is \textbf{visual} (a concept introduced by Schroeder - see \cite{BoS} or \cite{Buy})
 if there is a constant $D > 0$ (depending on the base point $a$ of $X$) such that for each $x\in X$ there is an infinite geodesic ray $\xi$ at $a$ such that
 $(x,\xi)_a > d(x,a)-D$.
\end{Def}

Recall $(x,\xi)_a=\liminf\limits_{n\to\infty}(x,\xi(n))_a$. Also recall
one can extend Gromov product $(x,y)_a$ to elements of $\partial X$ (see \cite[Sect.3]{BoS}):
For points $p,q\in \partial X$, the Gromov product is defined by
$$(p,q)_a= \inf \liminf\limits_{i\to\infty}(x_i,y_i)_a,$$
where the infimum is taken over all sequences $x_i,y_i$ such that $p=[x_i]$, $q=[y_i]$. Note that $(p,q)_a$ takes values in $[0, \infty]$,
and that ($p,q)_a=\infty$ if and only if $p=q$.

\begin{Prop}
Suppose  $f\colon X \to Y$ is a large scale Lipschitz function of proper hyperbolic spaces.
If $X$ is visual and there exist $\l_2,\m_2>0$ so that for every geodesic ray $\g\colon [0,\infty)\to X, \g(0)=a$ one has
$$
\l_2 \cdot d(x,y)-\m_2 \leq d(f(x),f(y)), \qquad \forall  x,y\in \im \g.
$$
then $f$ is radial.
\end{Prop}
\dokaz Suppose $f$ is $(\lambda,\mu)$-Lipschitz. Fix a constant $D > 0$ such that for each $u\in X$ there is an infinite geodesic $\xi$ such that
 $(u,\xi)_a > d(u,a)-D$, where $a$ is the base point of $X$. We claim existence of a universal (at $a$) constant $\mu_3$
 such that
 $$
\l_2 \cdot d(x,y)-\m_3 \leq d(f(x),f(y)), \qquad \forall  x,y\in au.
$$
There reason for it is existence on the ray $\xi$ of points $x_1,y_1$ such that $d(x,x_1) < D+\delta$ and $d(y,y_1) < D+\delta$.
Therefore $d(f(x),f(y))\ge d(f(x_1),f(y_1))-d(f(x),f(x_1))-d(f(y),f(y_1))\ge \l_2 \cdot d(x_1,y_1)-\m_2 -\lambda(2D+2\delta)-2\mu\ge
\l_2 \cdot d(x,y)-\l_2(2D+2\delta)-\m_2 -\lambda(2D+2\delta)-2\mu$
and $\mu_3=\l_2(2D+2\delta)+\m_2+\lambda(2D+2\delta)+2\mu$ works.
\edokaz

\section{Induced radial functions}
This section is devoted to the issue of H\" older maps between Gromov boundaries inducing radial functions between corresponding
hyperbolic spaces.

\begin{Def}
 A metric $d$ on $\partial X$ is \textbf{visual} (see \cite{BS} or \cite{BoS}) if there are constants $K, C > 1$ such that
$$\frac{K^{-(p,q)_a}}{C} \leq d(p,q)\leq C\cdot K^{-(p,q)_a}$$
Visual metrics exist for $K,C$ sufficiently large and they are all H\"older-equivalent.
\end{Def}

\begin{Lem} \label{CharOfHolderLemma}
Suppose $X$ and $Y$ are proper hyperbolic spaces with base points $a$ and $b$, respectively.
 A function $g:\partial X\to \partial Y$ is H\"older with respect to visual metrics
 if and only if there exist constants $A, B > 0$ such that for all geodesic rays $\xi_1,\xi_2$ in
 $X$ at $a$ one has
 $$(\eta_1,\eta_2)_b\ge A\cdot (\xi_1,\xi_2)_a-B$$
 for any geodesic rays $\eta_1,\eta_2$ at $b$ satisfying $g([\xi_i])=[\eta_i]$, $i=1,2$.
\end{Lem}
\dokaz
Choose visual metrics $d_X$ on $\partial X$ and $d_Y$ on $\partial Y$ with the same parameters $K,C$.
Assume $g$ is H\"older, hence there are constants $A,B$ such that $d_Y(g(p),g(q))\leq B\cdot d_X(p,q)^A$.
Therefore
$$\frac{K^{-(g(p),g(q))_b}}{C}\leq B\cdot (C\cdot K^{-(p,q)_a})^A=(B\cdot C^A)\cdot K^{-A(p,q)_a}$$
leading to existence of a constant $B_1$ such that
$$(g(p),g(q))_b\ge A\cdot (p,q)_a-B_1$$
for all $p,q\in\partial X$. Since Gromov products of elements of the boundary and Gromov products
of geodesic rays representing them differ at most by $\delta$, there exists a constant $B_2 > 0$ such that for all geodesic rays $\xi_1,\xi_2$ in
 $X$ at $a$ one has
 $$(\eta_1,\eta_2)_b\ge A\cdot (\xi_1,\xi_2)_a-B_2$$
 for any geodesic rays $\eta_1,\eta_2$ at $b$ satisfying $g([\xi_i])=[\eta_i]$, $i=1,2$.

Choose constants $A, B > 0$ such that for all geodesic rays $\xi_1,\xi_2$ in
 $X$ at $a$ one has
 $$(\eta_1,\eta_2)_b\ge A\cdot (\xi_1,\xi_2)_a-B$$
 for any geodesic rays $\eta_1,\eta_2$ at $b$ satisfying $g([\xi_i])=[\eta_i]$, $i=1,2$.

Notice $(g(p),g(q))_a \ge A\cdot (p,q)_a -B-(A+1)\delta$ for all $p,q\in \partial X$.
Therefore
$$d_Y(g(p),g(q))_a\leq  C\cdot K^{-(g(p),g(q))_a}\leq C\cdot K^{B+(A+1)\delta} \cdot K^{-A\cdot (p,q)_a}\leq M\cdot d_X(p,q)^A$$
for some constant $M$ independent on $p$ and $q$.
\edokaz

\begin{Def}
 Suppose $X$ and $Y$ are proper geodesic spaces with base points $a$ and $b$, respectively, and
 $g:\partial X\to \partial Y$ is a function. If $X$ is visual, then a \textbf{radial extension} $f:X\to Y$ of $g$ with parameter $A > 0$ is defined as follows:\\
 1. $f(a)=b$.\\
 2. If $x\in X\setminus{a}$ belongs to an infinite geodesic ray $\xi_x$ emanating at $a$,
 we pick $\eta_x$ emanating from $b$ and satisfying $[\eta_x]=g(\xi_x])$ and we put $f(x)=\eta_x(A\cdot t_x)$, where $x=\xi_x(t_x)$.\\
 3. For the remaining points of $X$ we pick a maximal geodesic ray $\zeta_x:[0,M_x]\to X$ at $a$ containing $x$,
 we pick an infinite geodesic ray $\xi_x$ at $a$ satisfying $(\zeta_x(M_x),\xi_x)_a > d(\zeta_x(M_x),a)-D$,
  we pick $\eta_x$ emanating from $b$ and satisfying $[\eta_x]=g([\xi_x])$,
 and we put $f(x)=\eta_x(A\cdot t_x)$, where $x=\zeta_x(t_x)$.
\end{Def}

\begin{Thm}\label{InducedRadialMaps}
 Suppose $X$ and $Y$ are proper Gromov hyperbolic spaces whose boundaries are equipped with visual metrics.
 If $g:\partial X\to \partial Y$ is a H\" older map and $X$ is visual, then there is $M > 0$ such that each
 radial extension  $f:X\to Y$ of $g$ with parameter $A\leq M$ is a radial function.
 \end{Thm}
\dokaz
Choose constants $M, B > 0$ such that for all geodesic rays $\xi_1,\xi_2$ in
 $X$ at $a$ one has
 $$(\eta_1,\eta_2)_b\ge M\cdot (\xi_1,\xi_2)_a-B$$
 for any geodesic rays $\eta_1,\eta_2$ at $b$ satisfying $g([\xi_i])=[\eta_i]$, $i=1,2$.
 Consider a radial extension $f$ of $g$ with parameter $A\leq M$. All we need to show is that $f$
is large scale Lipschitz.
Assume $r =d(x,y)$, $x,y$ lie on geodesic rays $\xi_1$ and $\xi_2$,
 respectively, and those geodesics are mapped by $f$ to geodesics $\eta_1$ and $\eta_2$ passing through
 $f(x)$ and $f(y)$, respectively. If the distance of one of $f(x),f(y)$ to $b$ is at most
 $\eta_1,\eta_2)_b$, then $d(f(x),f(y))\leq \delta+|d(f(x),b)-d(f(y),b)|\leq \delta+A\cdot r\leq Mr+\delta$.
 If $d(f(x),b),d(f(y),b)\ge (\eta_1,\eta_2)_b$,
 then
  $2(f(x),f(y))_b+2\delta\ge 2(\eta_1,\eta_2)_b\ge 2M\cdot (\xi_1,\xi_2)_a-2B\ge 2M\cdot (x,y)_a-2B\ge
 d(f(x),b)+d(f(y),b)-Mr-2B\ge 2(f(x),f(y))_b-d(f(x),f(y))-Mr-2B$
 leading to $d(f(x),f(y))\leq Mr+2B+2\delta$.
 \edokaz

\section{Induced H\"older map on the Gromov boundary}

In this section we prove that radial functions are visual and induce H\" older maps on geodesic
boundaries (when equipped with visual metrics).


Let us start by recalling the stability of geodesics as stated in \cite[Theorem 1.3.2]{BS}.

\begin{Thm} \label{ThmClassicStabilityOfGeodesics}[Stability of geodesics \cite[Theorem 1.3.2]{BS}]
Let $X$ be a $\d$-hyperbolic geodesic
space, $a \geq 1, b \geq 0$. There exists $H = H(a, b, \d) > 0$ such
that for every $N \in \NN$ the image of every $(a, b)$-quasi-isometric embedding
$f \colon \{1, \ldots, N\} \to X, im(f)$, lies in the $H-$neighborhood of any geodesic
$c\colon[0, l]\to X$ with $c(0) = f(1), c(l) = f(N)$, and vice versa, $c$ lies
in the $H-$neighborhood of $im(f)$.
\end{Thm}

A common application of Theorem \ref{ThmClassicStabilityOfGeodesics} states that given a quasi-isometry between hyperbolic geodesic spaces, there exists a constant $H$ so that the image of every geodesic of $X$ is at Hausdorff distance at most $H$ from the 'corresponding' geodesic in $Y$. In the case of radial functions we obtain the following, slightly weaker statement.

\begin{Thm}\label{ThmOurStabilityOfGeodesics}[Stability of pointed rays]
Suppose $f\colon X \to Y$ is a $(\l_1,\m_1)-$LSL $(\l_2, \m_2)$-radial function between $\d-$hyperbolic geodesic spaces. Then for every $x_0\in X$ there exists a constant $H=H(\l_1, \m_1, \l_2, \m_2, \d, x_0)$ so that the following holds: if $\g\colon[0,l]\to X, 0\mapsto x_0$ is a geodesic in $X$ and $\g'$ is a geodesic in $Y$ between $f(x_0)$ and $f(\g(l))$ then $d_H(\im f(\g),\g')\leq H$, where $d_H$ is the Hausdorff distance.
\end{Thm}

\dokaz
Observing that every $(\l_1,\m_1)-$LSL $(\l_2, \m_2)$-radial function $f \circ \g \colon[0,l]\to Y, 0\mapsto f(x_0)$ is a quasi-isometric embedding with parameters depending on the collection $\l_1, \m_1, \l_2, \m_2, x_0$ (see also Proposition \ref{PropIndependenceOfQRofBasepoint}), the proof follows by Theorem \ref{ThmClassicStabilityOfGeodesics} as the $(\l_1+\m_1)$-chain $\{f(N)\}$ is at Hausdorff distance at most $\l_1+\m_1$ from $\im f(\g)$.
\edokaz

The following lemma is essentially from \cite{V}.

\begin{Lem} \label{Lemma1} [Essentially Lemma 2.9 of \cite{V}]
Suppose $a,b,c$ are points in a geodesic metric space and $\a$ is a geodesic between $b$ and $c$. Then $(b,c)_a\leq d(a,\a)$.
\end{Lem}

\dokaz
For every $p\in \a$ we have $d(b,c)=d(b,p)+d(p,c)$ hence
$$
2d(a,p)\geq d(a,c)-d(c,p)+d(a,b)-d(b,p)=d(a,c)+d(a,b)-d(c,b)=2(b,c)_a.
$$
\edokaz

\begin{Prop}[The lower bound for the Gromov product as distorted by the radial function] \label{PropBoundGromovProductViaQR}
Suppose $X,Y$ are $\d-$hyperbolic geodesic spaces, and $f\colon X \to Y$ is $(\l_1,\m_1)-$LSL $(\l_2, \m_2)-$radial with respect to $a$. Then there exist constants $A,B>0$ dependent on $\l_1,\m_1,\l_2,\m_2, \d, a$ so that for every $b,c\in X$ we have
$$
(f(b),f(c))_{f(a)} \geq A \cdot (b,c)_a-B
$$
In particular, $f$ is visual.
\end{Prop}

\dokaz
The following notation is needed for the proof:
\begin{itemize}
  \item $b^*$ is a point on the geodesic $ac$  determined by the condition $d(b^*,a)=(b,c)_a$;
  \item $b'$ is a point on the geodesic $f(a)f(c)$  determined by the condition $d(b',f(a))=(f(b),f(c))_{f(a)}$;
  \item $b''=f(b^*)$
  \item $c^*$ is a point on the geodesic $ab$  determined by the condition $d(c^*,a)=(b,c)_a$;
  \item $c'$ is a point on the geodesic $f(a)f(b)$  determined by the condition $d(c',f(a))=(f(b),f(c))_{f(a)}$;
  \item $c''=f(c^*)$
\end{itemize}

By Theorem \ref{ThmOurStabilityOfGeodesics} there exists $H=H(\l_1, \m_1, \l_2, \m_2, \d, a)$ and $y\in f(a)f(c)$ so that $d(b'',y)\leq H$ and $d(c'',f(a)f(b))\leq H$.

Case 1: $y\in f(a)b'$.
$$
(f(b),f(c))_{f(a)} =d(f(a),b') \geq d(f(a),y)\geq d(f(a),f(b^*))-d(b'',y)\geq
$$
$$
\geq \l_2 d(a,b^*)-\m_2-H=\l_2 (b,c)_a-\m_2-H
$$

Case 2: $y\in f(c)b'$.
Note that $d(f(c),b')=(f(a),f(b))_{f(c)}$ as $b'$ divides the geodesic $ac$ into two parts, whose lengths are the corresponding Gromov products.
By Lemma \ref{Lemma1} we have
$$
(f(a),f(b))_{f(c)}\leq d(f(c),f(a)f(b))\leq
$$
$$
\leq d(f(c),y)+d(y,b'')+d(b'',c'')+d(c'',f(a)f(b))\leq
$$
$$
\leq d(f(c),y)+H+\l_1 \d + \m_1+H=
$$
$$
\leq (f(a),f(b))_{f(c)}-d(y,b')+H+\l_1 \d + \m_1+H
$$
hence
$$
d(y,b')\leq 2H+\l_1 \d + \m_1.
$$
Consequently, we obtain
$$
(f(a),f(b))_{f(c)} =d(f(a),b')\geq d(f(a),b'')-d(b'',y)-d(y,b')\geq
$$
$$
\geq d(f(a),f(b^*))- H - (2H+\l_1 \d + \m_1) \geq \l_2 d(a,b^*)-\m_2 - (3H+\l_1 \d + \m_1)=
$$
$$
=\l_2 (b,c)_a-(\m_2 + 3H+\l_1 \d + \m_1).
$$

Combining both cases we obtain constants $A,B$ of required dependency.
\edokaz

\begin{Thm} [Induced H\" older map on the boundary]
Suppose $X$ and $Y$ are proper geodesic hyperbolic spaces. If $f \colon X \to Y$ is radial, then $f$ induces a H\" older map between boundaries
      $$
      \tilde f\colon \di X \to \di Y
      $$
      when equipped with visual metrics.
\end{Thm}

\dokaz
Fix $a\in X$.
Choose constants $A$ and $B$ such that
$$
(f(b),f(c))_{f(a)} \geq A \cdot (b,c)_a-B
$$
for all $b,c\in X$. Use \ref{CharOfHolderLemma}.
\edokaz

\section{Surjectivity and coarse surjectivity}

This section is devoted to the interplay of subjectivity of the induced function by $f:X\to Y$ on the boundaries and coarse subjectivity of $f$.

\begin{Def}
A function $f\colon X \to Y$ is \textbf{coarsely surjective} if $f(X)$ is coarsely dense with a parameter $S$ in $Y$, i.e., $\forall y\in Y \exists x\in X: d(f(x),y)\leq S$.
\end{Def}

\begin{Lem}\label{LemmaEqualityOfCloseSequences}
Suppose $X$ is a proper hyperbolic geodesic space. Assume $(x_n)$ and $(y_n)$ are two sequences and $(x_n)$ represents a point in $\di X$. If there exists $R$ so that $d(x_n,y_n)\leq R, \forall n$ then $(y_n)$ represent the point $[(x_n)]$ in $\di _X$.
\end{Lem}

\dokaz
$$
2(x_i,y_j)_{x_0} = d(x_i,x_0)+d(y_j,x_0)-d(x_i,y_j)\geq
$$
$$
\geq d(x_i,x_0)+d(x_j,x_0)-R-d(x_i,x_j)-R=2(x_i,x_j)_{x_0}-2R
$$
Hence  $\liminf_{i,j\to \infty}(x_i,x_j)_{x_0}=\infty$ implies $\liminf_{i,j\to \infty}(x_i,y_j)_{x_0}=\infty$, which concludes the proof.
\edokaz

\begin{Thm} \label{ThmPersistenceSurjectivity}
[Persistence of surjectivity from large scale to small scale]
Suppose $f\colon X \to Y$ is a visual function of proper $\d-$hyperbolic geodesic spaces. If $f$ is coarsely surjective, then $\tilde f$ is surjective.
\end{Thm}
\dokaz Fix $S > 0$ such that $dist(y,f(X)) < S$ for all $y\in Y$. Given $[\{y_n\}_{n\ge 1}]\in \partial Y$ choose
$x_n\in X$ satisfying $d_Y(y_n,f(x_n)) < S$ for all $n\ge 1$. Notice $d_X(x_n,a)\to\infty$ as $n\to \infty$.
By switching to a subsequence we may achieve $(x_i,x_j)_a\to \infty$ as $i,j\to\infty$.
That implies $[\{x_n\}_{n\ge 1}]\in \partial X$ and $\tilde f([\{x_n\}_{n\ge 1}])=[\{y_n\}_{n\ge 1}]$.
\edokaz

\begin{Thm} \label{ThmReverseSurjectivity}
Suppose $f\colon X \to Y$ is radial and $X$, $Y$ are proper $\d-$hyperbolic geodesic spaces. If $\tilde f$ is surjective and $Y$ is visual, then $f$ is coarsely surjective.
\end{Thm}
\dokaz Fix $a\in X$ and put $b=f(a)$. $Y$ being visual means that the union $Y_0$ of all infinite rays stemming from $b$
is coarsely dense in $Y$. Use \ref{ThmOurStabilityOfGeodesics} to find $H > 0$ such that for each
infinite ray $\eta$ at $a$, the $H$-neighborhood of its image contains all rays equivalent to $\tilde f([\eta])$.
Therefore the $H$-neighborhood of $f(X)$ contains $Y_0$ and $f(X)$ is coarsely dense in $Y$.
\edokaz

\section{Inducing $n-to-1$ maps of boundaries}

Another example of interplay between small scale properties of $\tilde f$ and large scale
properties of $f$ is the cardinality of point-inverses of $\tilde f$. A function is called $n$-to-1 if its point-inverses
contain at most $n$ points. There is a large scale analog of it defined by Miyata-Virk in \cite{MV}.

\begin{Thm}\label{InducedIsnTo1Thm}
 Suppose $f:X\to Y$ is a visual function of proper hyperbolic spaces and $n\ge 1$.
 The induced map $\tilde f:\partial X\to \partial Y$ is $n$-to-1 if and only if for each $r > 0$ there is
 $B > 0$ such that for any sequence of $n+1$ points $x_1,\ldots ,x_{n+1}$ in $X$
 satisfying $(x_i,x_j)_a < r$ for $i\ne j$ there are two different indices $k,m\leq n+1$
 for which $(f(x_k),f(x_m))_b \leq B$.
\end{Thm}
\dokaz Suppose $\tilde f:\partial X\to \partial Y$ is $n$-to-1 and there is $r > 0$
such that for each $m\ge 1$ one has a sequence of $n+1$ points $x^m_1,\ldots ,x^m_{n+1}$ in $X$
satisfying $(x^m_i,x^m_j)_a < r$ and $(f(x^m_i),f(x^m_j))_b > m$ for $i\ne j$.
We may assume each sequence $\{x^m_i\}$ converges to $a_i\in \partial X$
in which case $a_i\ne a_j$ for $i\ne j$ but $\tilde f(a_i)=\tilde f(a_j)$ for all $i,j\leq n+1$,
a contradiction.

Suppose for each $r > 0$ there is
 $B > 0$ such that for any sequence of $n+1$ points $x_1,\ldots ,x_{n+1}$ in $X$
 satisfying $(x_i,x_j)_a < r$ for $i\ne j$ there are two different indices $k,m\leq n+1$
 for which $(f(x_k),f(x_m))_b \leq B$.
Suppose $\tilde f(a_i)=\tilde f(a_j)=z$ for some sequence of different points $a_1,\ldots,a_{n+1}$ of $\partial X$.
Choose $r > 0$ such that $U(a_i,r-\delta)$ are mutually disjoint for $i\ne j$.
Pick $x_i\in U(a_i,r-\delta)\cap f^{-1}(U(z,B+\delta))$ for $i\leq n+1$.
Notice $(x_i,x_j)_a < r$ as otherwise $x_j\in U(a_i,r-\delta)$.
Without loss of generality assume $(f(x_1),f(x_2))_b \leq B$.
However, $(f(x_1),z)_b > B+\delta$ and $(f(x_2),z)_b > B+\delta$
resulting in $(f(x_1),f(x_2))_b > B$, a contradiction.
\edokaz

\begin{Def} [Coarse $n-to-1$ function]
A coarse function $f\colon X \to Y$ is coarsely $n-to-1$ (originally named (B)${}_n$ in \cite{MV}) if for every $R$ there exists $S$ so that if $\diam(A)<R$ for $A\subset Y$ then $f^{-1}(A)=B_1\cup\ldots\cup B_n$ with $\diam(B_i)<S, \forall i$ (i.e., with $\{B_i\}$ being an $S-$bounded family).
\end{Def}

Theorem \ref{InducedIsnTo1Thm} suggests the following:

\begin{Prop}
A bornologous function $f\colon X \to Y$ is coarsely $n-to-1$ if and only if for every $R$ there exists $S$ such that if
$d_X(x_i,x_j) > S$ for all $i\ne j$, where $x_1,\ldots,x_{n+1}$ is a sequence in $X$, then there exist two different indices
$k,m\leq n+1$ satisfying $d_Y(f(x_k),f(x_m)) > R$.
\end{Prop}
\dokaz
Suppose $f$ has the property that for every $R$ there exists $S$ such that if
$d_X(x_i,x_j) > S$ for all $i\ne j$, where $x_1,\ldots,x_{n+1}$ is a sequence in $X$, then there exist two different indices
$k,m\leq n+1$ satisfying $d_Y(f(x_k),f(x_m)) > R$. For any $R > 0$ and any $y\in Y$ choose a maximal $S$-net $A$
 in the point-inverse $X_y$ of
the $R$-ball $B(y,R)$ (that means distances between different points of $A$ are bigger than $S$ and $d(x,A)\leq S$ for all $x\in X_y$).
The cardinality of $A$ is at most $n$ and $X_y\subset B(A,S)$. Therefore $f$ is coarse and is coarsely $n$-to-1.

Suppose $f$ is coarsely $n$-to-1 and for every $R$ choose $S$ so that if $\diam(A)<R$ for $A\subset Y$ then $f^{-1}(A)=B_1\cup\ldots\cup B_n$ with $\diam(B_i)<S, \forall i$. Assume $A=\{x_1,\ldots,x_{n+1}\}$ is a sequence in $X$ so that $d_X(x_i,x_j) > S$ for all $i\ne j$
and $d_Y(f(x_i),f(x_j)) < R$ for all $i\ne j$. Notice $\diam(A) < R$, hence $f^{-1}(A)=B_1\cup\ldots\cup B_n$ with $\diam(B_i)<S, \forall i$.
There exists two different indices $k,m$ so that $x_k,x_m\in B_i$ for some $i$. That means $d_X(x_k,x_m) < S$, a contradiction.
\edokaz

\begin{Thm} \label{ThmPersistenceNto1Coarse}
Suppose a function $f\colon X \to Y$ of proper $\d-$hyperbolic geodesic spaces is visual
and has the property that $f\circ \alpha$ is a coarse embedding for each coarse embedding $\alpha:[0,\infty)\to X$. If $f$ is coarsely $n-to-1$  then $\tilde f$ is $n-to-1$.
\end{Thm}

\dokaz
Suppose there are $n+1$ geodesic rays $\b_1, \ldots, \b_{n+1}$ based at $x_0$, representing different points in $\di X$, whose equivalent classes (i.e., the corresponding points on the boundary) are mapped to the class of the same geodesic $\g$.

Claim 1: $\exists M: d_{GH}(f(\b_i),\g)\leq M$.

Fix $i$. The class $[f(\b_i)]$ is represented by the sequence $(f(\b_i(n)))_{n\in \NN}$. By Theorem \ref{ThmOurStabilityOfGeodesics} (the stability of geodesics) there exists $H$ so that
$$
d_{GH}\big(f(x_0)f(\b_i(m)),f(\b_i([0,m]))\big)\leq H, \quad \forall m.
$$
As the geodesic $\g_i$ representing the sequence $(f(\b_i(n)))_{n\in \NN}$ in $\di Y$ is obtained from compact geodesics $f(x_0)f(\b_i(m))$ in the limit process (using the standard argument with an application of the Arzela-Ascoli Theorem) we obtain
$$
d_{GH}\big(\g_i,f(\b_i([0,\infty)))\big)\leq H.
$$
As $[\g_i]=[\g]$ (meaning that $d_{GH}(\g_i,\g)$ is finite) and the set of indices $i$ is finite, claim 1 follows.

Choose $S$ so that $n+1$ points in $X$ with mutual distances bigger than $S$ cannot be mapped to points with mutual distances at most $4L$. As geodesic rays $\b_1, \ldots, \b_{n+1}$ represent different points in $\di X$ the distance $d(\b_i(t),\b_j(t))$ goes to $\infty$ as $t \to \infty$ for $i\neq j$. Choose $L$ so that $d(\b_i(t),\b_j(t))>3S, \forall t\geq L-S, \forall i\neq j$. Observe that $f^{-1}(B(\g(L),2M))$ contains points $\b_i(L), \forall i$. On the other hand, these points are at least $3S-$disjoint, a contradiction. Hence $\tilde f$ is n-to-1.
\edokaz

\begin{Thm}[Large scale dimension raising theorem for hyperbolic spaces]\label{ThmMain}
Suppose $X$, $Y$ are proper $\d-$hyperbolic geodesic spaces and $f\colon X \to Y$ is a radial function. If $f$ is coarsely $(n+1)$-to-1 and coarsely surjective then
$$
\dim(\di Y)\leq \dim (\di X)+n.
$$

Moreover, if $X$ and $Y$ are hyperbolic groups then
$$
\asdim(Y)\leq \asdim (X)+n.
$$
\end{Thm}

\dokaz
By Theorems \ref{ThmPersistenceSurjectivity} and \ref{ThmPersistenceNto1Coarse} the induced map $\tilde f$ between boundaries is onto and $(n+1)$-to-1. As boundaries are compact metric, $\tilde f$ is closed. By the classical dimension raising theorem we have $\dim(\di Y)\leq \dim (\di X)+n$.

The second part follows by the Gromov conjecture as proved in \cite{BL}.
\edokaz

\section{An explicit example}

In this section we provide a specific example of a coarsely 2-to-1 radial function $f\colon X \to Y$ whose induced map raises the dimension of the Gromov boundary. The induced map in the Gromov boundary will be the standard 2-to-1 surjection from the Cantor set to the unit circle. For basic facts used here see \cite{KB}. Many simple technical details are ommitted.

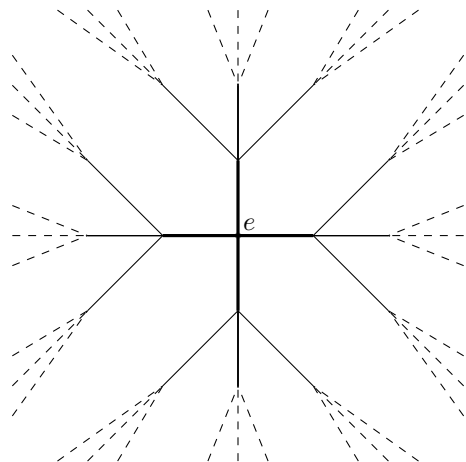
\begin{figure}
\begin{tikzpicture}[scale=.5]
\filldraw [gray] (0,0) circle (2pt);
\draw (0.3,0.3) node {$e$};
\draw[black, very thick] (-2,0)--(2,0);
\draw[black, very thick] (0,-2)--(0,2);
\foreach \x in {-1,1}
	{
	\foreach \y in {-1,0,1}
		{
		\draw (2 * \x,0)--(4*\x,2*\y);
		\draw (2 * \x,0)--(4*\x,0);
		\draw (0,2 * \x)--(2*\y,4*\x);
		\draw (0,2 * \x)--(0,4*\x);
		\foreach \z in {-1,0,1}
			{
			\draw [dashed](6 * \x,\z* 0.8 + 4*\y)--(4*\x,2*\y);
			\draw [dashed](\z* 0.8 + 4*\y, 6 * \x)--(2*\y, 4*\x);
			}
		}
	}
\end{tikzpicture}
\caption{The space $X$: the Cayley graph of the free group on two generators. Lengths of all edges equal $1$. }
\label{fig1}
\end{figure}

\begin{figure}
 \begin{tikzpicture}[scale=.5]
\filldraw [gray] (0,0) circle (2pt);
\draw (0.5,0.3) node {$x_0$};
\draw (0,0) circle (6);
\draw[black, very thick] (-2,0)--(2,0);
\draw[black, very thick] (0,-2)--(0,2);
\foreach \x in {1,..., 12}
	{
	\draw (\x * 360/12:2)--(\x * 360/12:3.5);
	}
\foreach \x in {1,..., 36}
	{
	\draw [dashed] (\x * 360/36:3.5)--(\x * 360/36:4.5);
	}
\end{tikzpicture}
\caption{The image of the map $f$ in the unit disc model of the hyperbolic plane.}
\label{fig2}
\end{figure}
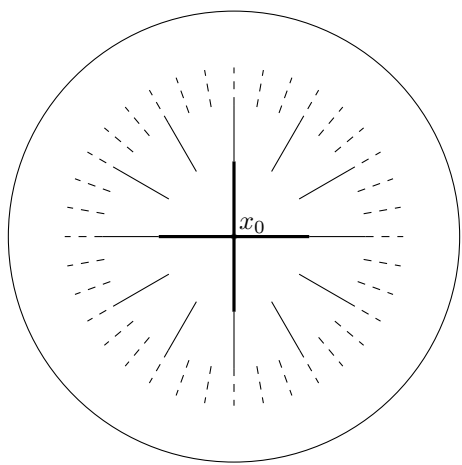

\textbf{Define $X$} (see Figure \ref{fig1}) to be the Cayley graph of the free group on two generators corresponding to the obvious representation with two generators and no relations. It is well known that its Gromov boundary is the Cantor set. In order to simplify the argument we will relabel the graph in the following way:
\begin{enumerate}
\item let vertex $e$ denote the base point of $X$;
\item vertices $V(X)$ are labeled by finite sequences $(a_i)$ of numbers from the set $\{0,1,2\}$ with $a_1$ possibly being $4$. Vertex $e$ corresponds to the empty sequence;
\item the vertex $e$ is of valence $4$, the corresponding edges (as well as their endpoints) are labeled by $0,1,2,3$;
\item inductively label the remaining edges: at each vertex labeled $w$ label the remaining three edges connecting $w$ to $(w,0), (w,1)$ and $(w,2)$ by $0,1$ and $2$ respectively. Notation $(w,i)$ denotes the sequence obtained from $w$ by adding en element $i$ at the end. Each vertex is of valence $4$.
\end{enumerate}

\textbf{Define $Y$} to be the Hyperbolic plane in the unit disc model. It is well known that its Gromov boundary is the unit circle. We will use the following notation:
\begin{enumerate}
\item $x_0$ is the centre of the disc;
\item for $r\geq 0$ and $\f\in\RR$ define  $A(r,\f)$ and $B(r,\f)$ to be points on the geodesic ray (a line) based at $x_0$ at angle $\f$ uniquely determined by conditions $d(A(r,\f),x_0)=r$ and $d(B(r,\f),x_0)=r+1$;
\item define $l(r,\f)$ to be a geodesic segment connecting $A(r,\f)$ to $B(r,\f).$
\end{enumerate}

\textbf{Define $f\colon X \to Y$} in the following way (see Figure \ref{fig2}):
\begin{enumerate}
\item $f(e)=x_0;$
\item for each $i$ map the edge $i$ emanating from $e$ isometrically to $l(0,i \cdot \pi/2)$;
\item suppose $f(w)=B(n,\f)$ for $w\in V(X)-\{x_0\}$ (in this case $n-1$ is the length of $w$).  For every $i\in \{0,1,2\}$ map the edge $i$ emanating from $w$  (excluding the origin $w$) isometrically onto $l(n,\f -\frac{2}{6} \frac{\pi}{2}3^{-(n-1)} + i \cdot \frac{2}{6} \frac{\pi}{2}3^{-(n-1)})$ (again, excluding $B(n,\f)$). Note that this causes a discontinuity of $f$ and that vertices $(w,i)$ are again mapped to points of form $B(n,\f_i)$;
\item in particular, $f$ maps edges of $X$ at distance $n$ from $e$ onto radial geodesic segments of length $1$ at  distance $n$ from $x_0 $, and distributes them uniformly with respect to angle.
\end{enumerate}

\textbf{The function $f$ is radial:} Suppose $x,y$ are points on a geodesic ray in $X$ based at $e$. Since $d(f(x),f(e))=d(x,e)$ and $d(f(y),f(e))=d(y,e)$ we have $d(f(x),f(y))\geq d(x,y)$, i.e., $f$ is $(1,0)-radial.$

\textbf{Discontinuities are uniformly bounded:} Roughly speaking, the discontinuities/jumps appear at distance $n$ from $x_0$ for $n\in\NN$. At fixed $n$ they are of the same size, that size being the shortest distance between two points of the form $A(n,\f)$; such points are uniformly distributed along the circle of radius $n$ around $x_0$. Thus the size of discontinuities is at most $c_n/N_n$, where $c_n$ is the circumference of the mentioned circle and $N_n$ is the number of discontinuities at radius $n$. It is well known that $c_n=2 \pi \sinh(n)$ and we can easily see that $N_n=4 \cdot 3^{n-1}$ meaning that $\lim_{n\to\infty}c_n/N_n<\infty$ which proves the claim.

\textbf{The function $f$ is LSL:} The function $f$ is an isometry on every edge, discontinuities on any given geodesic are at distance $1$ apart from each other. Since they are bounded by the previous paragraph by, say, $B'\in \RR$ we can easily deduct that $f$ is $(B'+1,B')$-LSL.

\textbf{The function $f$ is coarsely 2-to-1:}  proving this property is lengthy but easy and geometrically apparent. The reader is referred to a similar argument in \cite{MV}.

\textbf{Conclusion:} The function $f$ induces a 2-to-1 surjection between Gromov boundaries (from the Cantor set onto the circle) thus raising their dimension (and also the asymptotic dimension of the underlying space) by $1$.

\section{Inducing coarsely $n-to-1$ functions}
This section is devoted to general examples of $n$-to-1 maps of boundaries of hyperbolic spaces
inducing coarsely $n$-to-1 functions between hyperbolic spaces.

\begin{Def}
 Given a function $f:X\to Y$ of based hyperbolic spaces and $n\ge 1$ define
 $c_{f,n}:(0,\infty)\to (0,\infty]$ as follows:\\
 $c_{f,n}(r)$ is the supremum of numbers $B$ such that there is a sequence of $n+1$ points $x_1,\ldots ,x_{n+1}$ in $X$
 satisfying $(x_i,x_j)_a < r$ and $(f(x_i),f(x_j))_b \ge B$ and for all $i\ne j$.
\end{Def}

\begin{Thm}
Let $X$, $Y$ be proper hyperbolic spaces and $n\ge 1$
 Suppose $f:X\to Y$ is a radial extension with parameter $1$ of a H\" older function
between the boundaries.
 $f$ is coarsely $n$-to-1 if and only if
 there are numbers $r_0, M > 0$ such that $c_{f,n}(r)-r < M$ for all $r > r_0$.
 \end{Thm}
\dokaz
Suppose  $f$ is coarsely $n$-to-1and $c_{f,n}(r)-r$ is not bounded as $r\to \infty$. Let $S$ be the number
corresponding to $\delta+1$ ($f$ is coarsely $n$-to-1)
Find $r > 2S$ such that $c_{f,n}(r)-r > 2S+2$. Find a sequence of $n+1$ points $x_1,\ldots ,x_{n+1}$ in $X$
 satisfying $(x_i,x_j)_a < r$ and $(f(x_i),f(x_j))_b \ge r+2S+1$ and for all $i\ne j$.

Notice $d(x_i,a) > r+2S$ for each $i$ as otherwise $(f(x_i),f(x_j))_b \leq r+2S$ for any $j$.
Let $a_i$ be on a geodesic $[a,x_i]$ at the distance $r+2S$ from $a$.
$f(a_i)$ is on a geodesic $[b,f(x_i)]$ at the distance $r+2S$ from $a$. Therefore
$d(f(a_i),f(a_j))\leq \delta$ for $i\ne j$. As $f$ is coarsely $n$-to-1, there is a pair of indices
$k\ne m$ such that $d(a_k,a_m) < S$. Hence $(a_k,a_m)_a > r+2S-S/2$ contradicting
$(x_k,x_m)_a < r$ and $(a_k,a_m)_a\leq (x_k,x_m)_a$.

Suppose $c_{f,n}(r)-r < M$ for large $r$. Assume there are $n+1$ points $x_1,\ldots, x_{n+1}$
in $X$ with mutual distances at least $D$ being mapped to an $R$-ball $B(y,R)$ for some $y$ with $d(y,b) > K$.
Therefore $K-R < d(x_i,a) < K+R$ for all $i$ resulting in $(x_i,x_j)_a < K+R-D/2$. Also $(f(x_i),f(x_j))_b\ge K-R-R/2$ for all $i\ne j$.
That means $c_{f,n}(K+R-D) \ge K-2R$ and $c_{f,n}(K+R-D)-(K+R-D) > K-2R-K-R+D=D-3R$.
If $R$ is fixed, we can put $D=3R+M+1$ and
$K$ so large that $c_{f,n}(r)-r < M$ for $r \ge K+R-D$. That makes existence of points $x_1,\ldots, x_{n+1}$ above
impossible. Add the fact the point inverse of $B(a,K+2R)$ is bounded and conclude $f$ is coarsely $n$-to-1.
\edokaz

\begin{Thm}\label{FiniteTo1Theorem}
 Suppose $g:Z_1\to Z_2$ is a Lipschitz map of compact metric spaces and $n\ge 1$.
The induced radial function $f:Cone(Z_1)\to Cone(Z_2)$ between the corresponding hyperbolic cones
 is coarsely $n$-to-1 if there is $A > 0$ such that for any sequence of $n+1$ points $x_1,\ldots ,x_{n+1}$ in $Z_1$
 satisfying $d(x_i,x_j) > r$ for $i\ne j$ there are two different indices $k,m\leq n+1$
 for which $d(f(x_k),f(x_m)) \ge A\cdot r$.
\end{Thm}
\dokaz By rescaling the metric on $Z_2$ we may assume $g$ is $1$-Lipschitz.
As shown in \cite{Buy}, $Z_i$ is the Gromov boundary of $Cone(Z_i)$, $i=1,2$, and the metrics on $Z_1$ and $Z_2$ are visual with parameter $e$.
Hence, in terms of Gromov product the assumptions of the theorem translate as follows:
there are constants $B, C > 0$ such that $(g(z_1),g(z_2))_b\ge (z_1,z_2)_a-B$ for all $z_1,z_2\in Z_1$
and for each sequence of $n+1$ points $x_1,\ldots ,x_{n+1}$ in $Z_1$
 satisfying $(x_i,x_j)_a < r$ for $i\ne j$ one has $(g(x_k),g(x_m)) \leq r+C$ for some indices $k\ne m$.

 Consider the induced radial function
$f:Cone(Z_1)\to Cone(Z_2)$ with parameter $1$. Suppose there is a set of $n+1$ points $x_i\in Cone(Z_1)$
such that $(x_i,x_j)_a < r$ for some $r > 0$ and all $i\ne j$ and $(f(x_i),f(x_j))_b\ge r+C+\delta$ for all $i\ne j$.
Let $\xi_i$ be the geodesic ray from $a$ and passing through $x_i$ for $i\leq n+1$.
Since $(g(\xi_i),g(\xi_j))_b\ge (f(x_i),f(x_j))_b\ge r+C+\delta$ for all $i\ne j$, there is a pair $k,m$ of different
indices satisfying $(\xi_k,\xi_m)_a\ge r+\delta$. The only possibility for $(x_k,x_m)_a < r$
is if one of the points $x_k,x_m$ is at the distance at most $r+\delta$ from $a$.
In that case $(f(x_k),f(x_m))_b \leq r+\delta$ as well contradicting $(f(x_k),f(x_m))_b\ge r+C+\delta$.
\edokaz

\begin{Example}\label{Cantor2To1Example}
 The map $f:C\to I$, $f(\{a_i\})=\sum\limits_{i=1}^\infty \frac{a_i}{2^i}$, from the Cantor set represented as the product $\prod\limits_{i=1}^\infty \partial I$
 (with the metric $d(\{x_i\},\{y_i\})=\sum\limits_{i=1}^\infty \frac{|x_i-y_i|}{2^i}$)
 is $1$-Lipschitz and
has the property that for any sequence $x_1,x_2 ,x_3$ of $3$ points in $C$
 satisfying $d(x_i,x_j) > r$ for $i\ne j$ there are two different indices $k,m\leq 3$
 for which $d(f(x_k),f(x_m)) \ge r/2$.
\end{Example}
\dokaz $f$ being $1$-Lipschitz is obvious.
First consider $r=\frac{1}{2^m}$, $m\ge 1$. Given three points $x_1,x_2 ,x_3$ in $C$
satisfying $d(x_i,x_j) > r$ we may assume $x_1$ and $x_2$ form the closest pair. Let $k$ be the first index their coordinates
are different. Notice $k\leq m$. Therefore $x_3$ differs with both $x_1$ and $x_2$ on an earlier coordinate $p$.
Without loss of generality we may assume $x_3$ agrees with $x_1$ on $k$-th coordinate, so we might as well assume it is $0$.
Therefore $d(f(x_1),f(x_3))\ge \frac{1}{2^p}-\sum\limits_{i=p+1}^{k-1}\frac{1}{2^i}-\sum\limits_{i=k+1}^{\infty}\frac{1}{2^i}= \frac{1}{2^k}\ge \frac{1}{2^m}$.

For $r > 1$ the condition is vacuously true and for $r \leq 1$ we pick $m\ge 1$ such that
$\frac{1}{2^m} < r \leq \frac{1}{2^{m-1}}$. Given three points $x_1,x_2 ,x_3$ in $C$
satisfying $d(x_i,x_j) > r$ we know there is a pair satisfying $d(f(x_k),f(x_m)) \ge \frac{1}{2^m}> r/2$.
\edokaz

\begin{Cor}
 There is a radial coarsely $2$-to-1 function $f:X\to Y$ of proper hyperbolic spaces
 such that $asdim(X)=1$ and $asdim(Y)=2$.
\end{Cor}
\dokaz Consider the $2$-to-1 map $f$ from Example \ref{Cantor2To1Example} and the induced radial
function $g$ with parameter $1$ of the hyperbolic cone $X=Cone(C)$ to the hyperbolic cone $Y=Cone(I)$.
$asdim(X)\ge 1$ as $X$ contains infinite rays. To assert $asdim(X)\leq 1$ we invoke the theorem of S.Buyalo \cite{Buy}:
$asdim(X)$ is at most capacity dimension of $\partial X$ plus 1. Recall (see Proposition 3.2 of \cite{Buy}) that the capacity dimension of a metric space $Z$ is the
infimum of all integers $m$ with the following property: There exists a constant $c > 0$ such that for all sufficiently small $s > 0$, $Z$ has a $cs$-bounded covering with $s$-multiplicity at most $m + 1$ (meaning that an $s$-ball intersects at most $m+1$ elements of the covering). The product $\prod\limits_{i=1}^\infty \partial I$ has covers of the form $const \times \prod\limits_{i=n}^\infty \partial I$
($n$ fixed) that are $4\cdot 2^{-n-1}$-bounded and have $2^{-n-1}$-multiplicity at most $1$.

Using the same theorem one gets $asdim(Y)\leq 2$. The reason for $asdim(Y)\ge 2$ is that $Y$ contains topological $2$-disks approaching
infinity.
\edokaz

\end{document}